\documentclass[12pt]{article}

\usepackage[latin1]{inputenc}
\usepackage[english]{babel}
\usepackage{amsmath,amsfonts,amssymb}
\usepackage{graphicx}
\usepackage{epstopdf}
\usepackage{xcolor}
\usepackage{tikz}
\usetikzlibrary{shapes}
\usetikzlibrary{arrows.meta}
\tikzset{%
  base/.style = {rectangle, rounded corners, draw=black,
    minimum width=1cm, minimum height=0.5cm,
    text centered, font=\sffamily},
  process/.style = {base, minimum width=0.5cm, fill=orange!15, font=\ttfamily},
  process1/.style = {rectangle, rounded corners, draw=black}
}
\usepackage[]{algorithm2e}
\usepackage{float}
\usepackage{eurosym}

\usepackage{booktabs}
\usepackage{hyperref}

\usepackage{commandes}

\def\x{\boldsymbol{X}}
\def\u{\boldsymbol{U}}
\def\alea{\boldsymbol{W}}
\def\w{\boldsymbol{W}}
\def\fw{\overline{w}}

\def\next{t+1}
\def\post{{t+1}}

\def\final{T}

\def\Demandel{\mathbf{D}^{el}}

\def\Ti{\boldsymbol{\theta}^i}
\def\ti{\theta^i}
\def\Tw{\boldsymbol{\theta}^w}

\def\tank{H}
\def\price{\pi}

\def\opt{^{\sharp}}

\def\price{p}
\def\transfert{\mathbf{F}}
\def\Demandth{\mathbf{D}^{th}}

\def\PV{\boldsymbol{\Phi}^{pv}}

\def\uad{\cU^{ad}}
\def\tw{\theta^w}
\def\to{\theta^e}
\def\pint{\Phi^{int}}
\def\pext{\Phi^{ext}}

\def\tank{h}

\newcommand{\ClosedIntervalOpen}[2]{[#1,#2)}

\graphicspath{{img/},{img/devices}}
\newcommand{\keywords}[1]{}
\usepackage{fullpage}

\title{Optimization of a domestic microgrid equipped with solar panel and battery: Model Predictive Control and Stochastic Dual Dynamic Programming approaches}

\author{Fran\c{c}ois Pacaud\footnote{CERMICS, Ecole des Ponts, Marne-la-Vall\'ee, France},
  \and
  Pierre Carpentier\footnote{UMA, ENSTA Paris, Institut Polytechnique de Paris, Palaiseau, France},
  \and
  Jean-Philippe Chancelier\footnotemark[1],
  \and
  Michel De Lara\footnotemark[1]
}

\begin{document}
\maketitle

\begin{abstract}
  In this study, a microgrid with storage (battery, hot water tank) and solar
  panel is considered. We benchmark two algorithms, MPC and SDDP,
  that yield online policies to manage the microgrid,
  and compare them with a rule based policy.
  Model Predictive Control (MPC) is a well-known algorithm
  which models the future uncertainties with a deterministic forecast.
  By contrast, Stochastic Dual Dynamic Programming (SDDP) models the future
  uncertainties as stagewise independent random variables with known probability distributions.
  We present a scheme, based on out-of-sample validation, to fairly compare the
  two online policies yielded by MPC and SDDP. Our numerical studies put to
  light that MPC and SDDP achieve significant gains compared to the rule based policy,
  and that SDDP overperforms MPC not only on average but on most of the
  out-of-sample assessment scenarios.

\end{abstract}


\begin{table}
 \centering
 \begin{tabular}{rl}
   \hline
   \hline
   \multicolumn{2}{c}{Abbreviations} \\
   \hline
   \hline
   AR & Auto-regressive \\
   EMS & Energy Management System \\
   MPC & Model Predictive Control \\
   SDP & Stochastic Dynamic Programming \\
   SDDP & Stochastic Dual Dynamic Programming \\
   \hline
   \hline
   \multicolumn{2}{c}{Physical variables} \\
   \hline
   \hline
   $t$ & Time \\
   $\Delta$ & Time step (15mn) \\
   $T_0$ & Horizon (24h) \\
   $p_t^e$ & Electricity price (Euro \euro)\\
   $p_t^d$ & Thermal comfort (virtual price) \\
   $A_p$ & Surface of solar panel ($m^2$) \\
   $\rho_c$ & Battery charge rate \\
   $\rho_d$ & Battery discharge rate \\
   $\beta_\tank$ & Conversion yield for hot water tank \\
   \hline
   \hline
   \multicolumn{2}{c}{Decision variables and uncertainties} \\
   \hline
   \hline
   $(\Omega, \cF, \PP)$ & Probability space \\
   $\va b_t$ & Level of energy in battery (kWh) \\
   $\va h_t$ & Level of energy in hot water tank (kWh) \\
   $\Demandel_t$ & Electrical demand (kW) \\
   $\Demandth_t$ & Hot water demand (kW) \\
   $\PV_t$ & Production of the solar panel (kW) \\
   $\theta_t^i$ & Inner temperature (${}^\circ C$)\\
   $\theta_t^w$ & Walls temperature (${}^\circ C$) \\
   $\theta_t^e$ & Outdoor temperature (${}^\circ C$) \\
   $\va f_t^b$ & Energy exchanged with the battery (kW) \\
   $\va f_t^{w}$ & Energy injected in the hot water tank (kW) \\
   $\va f_t^{h}$ & Energy injected in the electrical heater (kW) \\
   $\va u_t= (\va f_t^b, \va f_t^{w}, \va f_t^{h})$ & Controls \\
   $\va w_t= (\Demandel_t, \Demandth_t, \PV_t)$ & Uncertainties \\
   $\va x_t= (\va b_t, \va h_t, \Tw_t, \Ti_t)$ & States \\
   \hline
   \hline
   \multicolumn{2}{c}{Mappings} \\
   \hline
   \hline
   $f_t: \XX_t \times \UU_t \times \WW_\post \rightarrow \XX_\post$ & Linear dynamics \\
   $L_t: \XX_t \times \UU_t \times \WW_\post \rightarrow \RR$ & Convex operational cost \\
   $K: \XX_\final \rightarrow \RR$ & Convex final cost \\
   $\pi_t: \XX_t \rightarrow \UU_t$ & Control policy \\
   \hline
   \hline
 \end{tabular}
 \caption{Nomenclature}
 \label{tab:nomenclature}
\end{table}

\section{Introduction}

\subsection{Background introduction}
\label{sec:intro:background}

A microgrid is a local energy network that produces part of its energy
and controls its own demand. Such systems are complex to control because,
on the one hand, of the different stocks and interconnections,
and, on the other hand,
of electrical demands and weather conditions
(heat demand and renewable energy production) that are highly variable and
hard to predict at local scale
(see~\cite{olivares2014trends,morstyn2016control} for a panorama
of the challenges faced when controlling microgrids).

We consider here a domestic microgrid
equipped with a battery, an electrical hot water tank and a solar panel,
as in Figure~\ref{fig:twostocksschema}.
The microgrid is connected to an external grid to import electricity when
needed.
The battery stores energy when external grid prices are low
or when the production of the solar panel is above the electrical demand.
The house's envelope also plays the role of heat storage.
As a consequence, the system has four stocks to store energy: a battery, a hot water tank,
and two passive stocks being the house's walls and inner rooms.
Two kinds of uncertainties affect the system:
the electrical and domestic hot water demands are not known in advance;
the production of the solar panel is substantially perturbed by the variable weather nebulosity.

We aim to compare two classes of algorithms to tackle uncertainties
in a microgrid Energy Management System (EMS).
The Model Predictive Control (MPC) algorithm
(or its stochastic variant, Stochastic Model Predictive Control)
relies on a mathematical representation of the future uncertainties
as deterministic forecasts; then, MPC computes decisions online
as solutions of a deterministic multistage optimization problem.
Stochastic Dual Dynamic Programming (SDDP)
relies on a mathematical representation of the future uncertainties
as stagewise independent random variables with known probability distributions;
then, SDDP computes offline a set of value functions by backward induction,
and computes online decisions as solutions of a single stage stochastic
optimization problem, using the value functions. We
present a fair comparison of these two algorithms, and highlight the
pros and cons of both methods.

\subsection{Literature review}
\label{sec:intro:litterature}

\subsubsection*{Optimization and energy management systems}
EMS are integrated automated tools used to
monitor and control energy systems. 
In~\cite{olivares2014trends}, the authors
give an overview of the use of optimization methods in the design of EMS.
The MPC algorithm~\cite{garcia1989model} and its stochastic
variant, Stochastic Model Predictive Control (SMPC)~\cite{mesbah2016stochastic},
have been widely used to control EMS.
We refer the reader
to~\cite{vsiroky2011experimental}
for applications of MPC in buildings.
In~\cite{holjevac2015adaptive}, the authors use MPC for the optimal control
of a domestic microgrid, and investigate how to balance  the
uncertainties of renewable energies with the microgrid.
In~\cite{olivares2015stochastic}, the authors apply SMPC
to the management of an isolated microgrid, and highlight the benefit
of this method; an application of SMPC in buildings is presented
in~\cite{appino2018use};
a variant based on robust optimization is proposed in~\cite{paridari2016robust}.

\subsubsection*{Stochastic optimization}
At local scale, electrical demand and production are highly
variable, especially as microgrids are expected to absorb renewable energies.
This leads to pay attention to stochastic optimization
approaches~\cite{delara2014cfe}. Stochastic optimization has been widely applied to hydrovalleys
management~\cite{pereira1991multi}. Other applications have arisen recently,
such as integration of wind energy and
storage~\cite{haessig2015energy} or insulated microgrids
management~\cite{heymann2016stochastic,alasseur2019regression}.

Stochastic Dynamic Programming (SDP)~\cite{bertsekas1995dynamic} is a general
method to solve stochastic optimal control problems.
In energy applications, a variant of SDP, Stochastic Dual Dynamic Programming
(SDDP), has demonstrated its adequacy for large scale {convex} applications. SDDP was first
described in the seminal paper~\cite{pereira1991multi};
we refer
to~\cite{shapiro2011analysis} for a generic description of the algorithm and
its application to the management of hydrovalleys;
a proof of convergence in the linear case is given
in~\cite{philpott2008convergence}, and in the convex case
in~\cite{girardeau2014convergence}.
Recent articles have applied SDDP to the management of energy systems.
In~\cite{bhattacharya2018managing},
numerical experiments show that SDDP yields better results than a
myopic policy.
In~\cite{papavasiliou2018application}, SDDP is applied to the dispatch
of energy inside the German national grid, under time correlated uncertainties;
the authors observe that SDDP achieves better performances than
those of a deterministic-based policy.
{
  Other Approximate Dynamic Programming algorithms have been designed to tackle
  different stochastic optimal control problems like, for instance,
  incorporating probability constraints~\cite{balata2021statistical}.
}

\subsection{Main contributions and structure of the paper}
\label{sec:intro:contributions}

We provide a rigorous mathematical formulation of the optimal management of a domestic microgrid
--- equipped with a battery, an electrical hot water tank and a solar panel,
and connected to an external supply network ---
under stochasticity of demand and of renewable energy production.
To manage the microgrid,
we design online policies using two different algorithms,
MPC and SDDP.
Then, we develop a fair benchmark methodology to compare the two algorithms,
based on a realistic use case.
The comparison reveals that SDDP overperforms MPC not only on average but,
interestingly, on most of the out-of-sample assessment scenarios.

The paper is organized as follows.
In Sect.~\ref{sec:optimization},
we detail the modeling of a small residential microgrid and formulate
a mathematical multistage stochastic optimization problem.
Then, we outline the two algorithms, MPC and SDDP, in Sect.~\ref{sec:algo}.
Finally, in Sect.~\ref{sec:numeric} we detail the benchmark methodology
 and we provide numerical results on
the systematic comparison of MPC, SDDP and a rule based policy.
Sect.~\ref{Conclusion} concludes and the Appendix~\ref{sec:problem} provides
details on the physical equations of the microgrid.

\section{Optimization problem statement}
\label{sec:optimization}

\begin{figure}[!ht]
  \begin{center}
    \begin{tikzpicture}[scale=1]
      \path (0,1.5) node[draw,shape=circle,fill=black] (elec) {};
      \path (0,-1.5) node[draw,shape=circle,fill=black] (therm) {};

      \path (2.5,1.5) node[draw,shape=circle,fill=black] (p2) {};
      \path (-2.5,1.5) node[draw,shape=circle,fill=black] (p3) {};
      \path (0,3.7) node[draw,shape=circle,fill=black] (p4) {};

      \path (2.5,3.7) node[draw,shape=circle,fill=black] (p7) {};

      \path (2.5,-1.5) node[draw,shape=circle,fill=black] (p5) {};
      \path (-2.5,-1.5) node[draw,shape=circle,fill=black] (p6) {};


      \draw[very thick, shorten <=-1pt,transform canvas={yshift=-0.1cm}, node distance=2cm, ->]
      (elec) -+ (0, -.5cm) node [pos=0.5, right] {$f^h$};

      \draw[very thick, shorten <=-1pt,transform canvas={xshift=-0.1cm}, node distance=2cm, <->]
      (elec) -+ (-1.7, 1.5) node [pos=0.5, below] {$f^b$};
      \draw[very thick, shorten <=-1pt,transform canvas={xshift=0.1cm}, node distance=2cm, ->]
      (elec) -+ (1.7, 1.5) node [pos=0.5, below] {$d^{el}$};
      \draw[very thick, shorten <=-1pt,transform canvas={yshift=0.1cm}, node distance=2cm, <-]
      (elec) -+ (0, 2.8) node [pos=0.5, left] {$f^{ne}$};

      \draw[very thick, shorten <=-1pt,transform canvas={xshift=0.5cm}, node distance=2cm, ->]
      (therm) -+ (1.35, -1.5) node [pos=0.5, below] {$d^{hw}$};
      \draw[very thick, shorten <=-1pt,transform canvas={xshift=-0.1cm}, node distance=2cm, ->]
      (elec) -+ (-1.8, -0.8) node [pos=0.5, below] {$f^t$};

      \draw[very thick, shorten <=-1pt,transform canvas={xshift=0.1cm}, node distance=2cm, <-]
      (elec) -+ (1.7, 3.) node [pos=0.5, above] {$\phi^{pv}$};

      \node[inner sep=0pt] (bulb) at (p2)
      {\includegraphics[width=.07\textwidth]{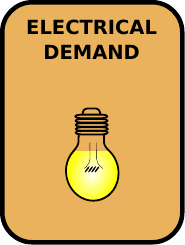}};
      \node[inner sep=0pt] (grid) at (p4)
      {\includegraphics[width=.07\textwidth]{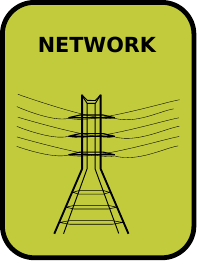}};
      \node[inner sep=0pt] (battery) at (p3)
      {\includegraphics[width=.07\textwidth]{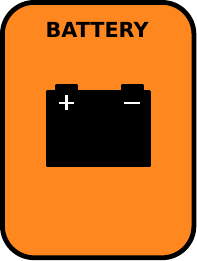}};
      \node[inner sep=0pt] (boiler) at (p7)
      {\includegraphics[width=.07\textwidth]{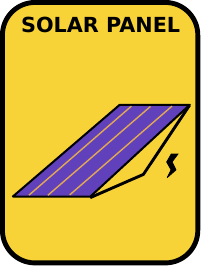}};

      \node[inner sep=0pt] (tank) at (therm)
      {\includegraphics[width=.07\textwidth]{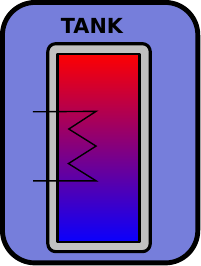}};
      \node[inner sep=0pt] (boiler) at (p6)
      {\includegraphics[width=.07\textwidth]{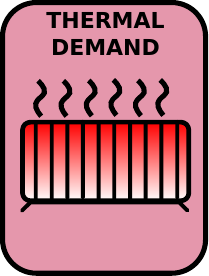}};
      \node[inner sep=0pt] (boiler) at (p5)
      {\includegraphics[width=.07\textwidth]{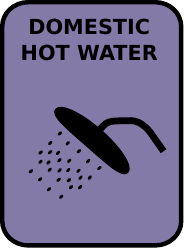}};

    \end{tikzpicture}
  \end{center}
  \caption{Electrical microgrid}
  \label{fig:twostocksschema}
\end{figure}
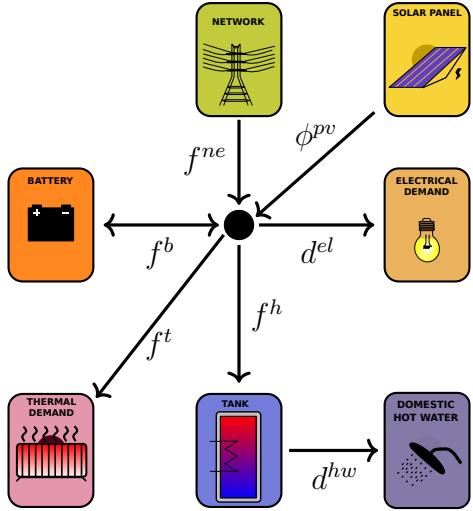

We consider the optimal management of a microgrid which,
as depicted in~Figure~\ref{fig:twostocksschema}, consists of a single house
equipped with a battery and an electrical hot water tank. An electrical heater
can produce heat in winter, and a solar panel can produce energy locally.
The decision maker aims at minimizing the energy bill ---
that is, the cost of the possible recourse energy supplied by the external network --- while
satisfying the energy demands (hot water and electricity),
and ensuring a minimal thermal comfort.

In this section, we write up a multistage stochastic optimization problem.
As decisions are taken at discrete time steps, we start by discretizing the
time interval in~\S\ref{Decisions_are_taken_at_discrete_times}. Then,
we introduce the uncertainties in~\S\ref{sec:model:alea},
the controls and the stocks in~\S\ref{sec:model:control}
and~\S\ref{sec:model:states}.
We detail the nonanticipativity constraints in~\S\ref{sec:model:cons1},
the bounds constraints in~\S\ref{sec:model:cons2} and
the objective function in~\S\ref{sec:objective}.
Finally, we formulate a multistage stochastic optimization problem in~\S\ref{sec:model:soc}.

\subsection{Decisions are taken at discrete time steps}
\label{Decisions_are_taken_at_discrete_times}

The EMS takes decisions every 15 minutes to control the system.
We consider a time interval~$\Delta = 15 \text{mn}$, a time horizon~$\final_0{=}24\text{h}$ and a
number~\( \final{=}\final_0/\Delta{=}96\) of time steps.
We adopt the following convention for discrete time processes:
for each time step~$t \in \{0, 1, \cdots, \final \}$,
$x_t$ denotes the value of the variable~$x$ at the
beginning of the interval~$\ClosedIntervalOpen{t \Delta}{(t+1) \Delta}$.
Otherwise stated, we will denote by $\ClosedIntervalOpen{t}{t+1}$ the continuous time interval
$\ClosedIntervalOpen{t \Delta}{(t+1) \Delta}$.\footnote{Here, by using the notation $\ClosedIntervalOpen{t}{t+1}$,
we mean \emph{the time interval between~$t$ and $t+1$, excluding
$t+1$}. Indeed, we denote a time interval
between two decisions by~\( \ClosedIntervalOpen{t}{t+1} \), and not
by~\( [t,t{+}1] \), to indicate that a decision is taken
at the beginning of the time interval~\( \ClosedIntervalOpen{t}{t+1} \),
and that a new one will be taken
at the beginning of the time interval~\( \ClosedIntervalOpen{t{+}1}{t{+}2} \),
and that these two consecutive intervals do not overlap.}

\subsection{Modeling uncertainties as random variables}
\label{sec:model:alea}

Because of their unpredictable nature, the EMS cannot anticipate the values of
electrical and thermal demands, nor the production of the solar panel.  We
choose to model these quantities as random variables over a probability space
$(\Omega, \cF, \PP)$. We adopt the following convention: a random variable will
be denoted by an uppercase bold letter~$\va Z$, and its realization
for a given outcome~$\omega\in \Omega$ will be
denoted in lowercase~$z = \va Z(\omega)$.  We denote by $\Demandel_t$ the
electrical demand, $\Demandth_t$ the hot water demand and $\PV_t$ the production
of the solar panel, all being real-valued random variables.  For each time
step~$t = 1, \ldots, \final$, we define the (uncertainty) random variable
\begin{equation}
  \w_t = (\Demandel_t, \Demandth_t, \PV_t)
  \eqfinp
  \label{uncertainties}
\end{equation}
The uncertainty $\w_t$ is a multivariate random variable taking values in~$\WW_t = \RR^3$.

\subsection{Modeling controls as random variables}
\label{sec:model:control}
As decisions depend on the previous uncertainties, controls are random
variables.
At the beginning of the time interval~\( \ClosedIntervalOpen{t}{t+1} \),
the EMS takes three decisions:
\begin{itemize}
  \item $\va f_t^b$, how much energy to charge in/discharge from the battery,
  \item $\va f_t^{w}$, how much energy to store in the electrical hot water tank,
  \item $\va f_t^{h}$, how much energy to inject in the electrical heater,
\end{itemize}
during the time interval~\( \ClosedIntervalOpen{t}{t+1} \).
For each time step~$t = 1, \ldots, \final$,
we define the decision multivariate random variable,
taking values in~$\UU_t =\RR^3$, as
\begin{equation}
  \u_t = (\va f_t^b, \va f_t^{w}, \va f_t^{h})
\eqfinp
\end{equation}
During the time interval~\( \ClosedIntervalOpen{t}{t+1} \), the EMS imports an energy
quantity~$\va f^{ne}_\post$ from the external network in order to fulfill the load balance
equation
\begin{equation}
  \label{eq:loadbalancerandom}
  \va f_\post^{ne}+\PV_\post =
\va f_t^b + \va f^{w}_t + \va f^{h}_t + \Demandel_\post
\eqfinv
\end{equation}
whatever the demand~$\Demandel_\post$ and the production of the solar panel
$\PV_\post$, unknown
at the beginning of the time interval~\( \ClosedIntervalOpen{t}{t+1} \).
On the left-hand side of Equation~\eqref{eq:loadbalancerandom}, the load consists of
\begin{itemize}
\item $\va f^{ne}_\post$, the energy surplus or shortage
    (when $\va f^{ne}_\post <0$, one wastes the surplus;
    when $\va f^{ne}_\post >0$, one imports the shortage from the network),
      \item $\PV_\post$, the production of the solar panel,
\end{itemize}
all during the time interval~\( \ClosedIntervalOpen{t}{t+1} \).
On the right-hand side of Equation~\eqref{eq:loadbalancerandom}, the electrical demand
is the sum of
\begin{itemize}
  \item $\va f^b_t$, the energy exchanged with the battery,
  \item $\va f^{w}_t$, the energy injected into the electrical hot water tank,
  \item $\va f^{h}_t$, the energy injected in the electrical heater,
  \item $\va d^{el}_\post$, the electrical demands (lightning, cooking...), aggregated in a single
    demand,
\end{itemize}
all during the time interval~\( \ClosedIntervalOpen{t}{t+1} \).

Later, we will aggregate the solar panel production~$\PV_\post$ with
the demands~$\Demandel_\post$ in Equation~\eqref{eq:loadbalancerandom},
as these two quantities appear only by their difference.


\subsection{States and dynamics\label{sec:state}}
\label{sec:model:states}

The state is the multivariate random variable
\begin{equation}
  \x_t = (\va b_t, \va h_t, \Tw_t, \Ti_t)
  \eqfinv
 \label{eq:state}
\end{equation}
which consists of the stocks~$\va b_t$ in the battery and~$\va h_t$ in the
electrical hot water tank, plus the two temperatures~$( \Tw_t, \Ti_t)$ of the
thermal envelope. Thus, the state random variable~$\x_t$ takes values in~$\XX_t = \RR^4$.

The discrete time dynamics
\( f_t : \XX_t \times \UU_t \times \WW_\post \rightarrow \XX_\post \)
describes the time evolution
\begin{equation}
\label{eq:dynamics}
    \x_{t+1} = f_t \bp{\x_t, \u_t, \alea_{t+1}}
\end{equation}
of the state,
where~$f_t$ is a piecewise linear function
(a property that will prove important for the SDDP algorithm),
that corresponds to the integration of the continuous
dynamics~\eqref{eq:batteryequation}-\eqref{eq:tankequation}-\eqref{eq:r6c2continuous}
given in Appendix~\ref{sec:problem}.
We suppose that we start from a given state~$x_0 \in \XX_0 = \RR^4$, thus adding
the constraint~$\x_0 = x_0$.

\subsection{Nonanticipativity constraints}
\label{sec:model:cons1}

The future realizations of uncertainties are unknown.
Thus, decisions at time step~$t$ are functions of previous history only, that is,
the information collected between step~$0$ and step~$t$. Such a constraint
is encoded as an algebraic constraint, using the tools of Probability
theory~\cite{kallenberg2002}, in the so-called
\emph{nonanticipativity constraints} written as
\begin{subequations}
\begin{equation}
  \label{eq:measurability}
  \sigma(\u_t) \subset \cF_t
\eqfinv
\end{equation}
where~$\sigma(\u_t)$ is the~$\sigma$-algebra generated by the random
variable~$\u_t$ and
$\cF_t = \sigma(\w_1, \cdots, \w_t)$ is the~$\sigma$-algebra
generated by the previous uncertainties~$(\w_1,$ $\ldots, \w_t)$.
If Constraint~\eqref{eq:measurability} holds true, the Doob Lemma~\cite{kallenberg2002}
ensures that there exists a measurable function
$\pi_t : \WW_1 \times \dots \times \WW_t \rightarrow  \UU_t$
such that
\begin{equation}
  \u_t = \pi_t(\w_1, \ldots, \w_t)
\eqfinp
  \label{eq:Doob}
\end{equation}
\end{subequations}
This is how we turn an (abstract) algebraic constraint into a more practical
functional constraint. The function~$\pi_t$ is called a \emph{policy}.

\subsection{Bounds constraints}
\label{sec:model:cons2}

The stocks~$\va b_t$ in the battery and~$\va h_t$ in the tank are bounded:
\begin{equation}
  \label{eq:statebounds}
  \underline b \leq \va b_t \leq \overline b \eqsepv
  0 \leq \va h_t \leq \overline h
\eqfinp
\end{equation}
At time step~$t$, the control~$\va f^b_t$ must ensure that the next state~$\va b_\post$
is admissible, that is, satisfies $\underline b \leq \va b_\post \leq \overline b$,
which, by time discretization of~\eqref{eq:batteryequation}, is equivalent to
the two inequalities\footnote{%
We have used the notation
$f^{+} = \max\na{0, f}$ and $f^{-} = \max\na{0, -f}$.}
\begin{equation}
  \underline b \leq
  \va b_t + \Delta \Bc{\rho_c (\va f_t^b)^+ + \dfrac{1}{\rho_d} (\va f_t^b)^-}
  \leq \overline b
  \eqfinv
\end{equation}
with~$\rho_c$ and~$\rho_d$ being the charge and discharge efficiencies of the battery.
Thus, the constraints on~$\va f_t^b$ depends on the stock~$\va b_t$.
The same reasoning applies for the tank energy~$\va f_t^{w}$ and the stock $\va h_t$.
Furthermore, we set bound constraints on controls:
\begin{equation}
  -\overline f^b \leq \va f_t^b \leq \overline f^b \eqsepv
  0 \leq \va f_t^{w} \leq \overline f_t^{w} \eqsepv
  0 \leq \va f_t^{h} \leq \overline f_t^{h} \eqfinp
\end{equation}
Finally, the load-balance equation~\eqref{eq:loadbalancerandom} also acts as a constraint
on the controls.
We gather all these constraints into an admissible subset
\( \uad_t(\x_t) \) of \(  \UU_t = \RR^3 \),
depending on the current state~$\x_t$, giving the constraint\footnote{%
  More formally, for each time step~$t = 1, \ldots, \final$,
      $\uad_t$ is a nonempty set-valued mapping $\uad_t: \XX_t \rightrightarrows \UU_t$.}

\begin{equation}
  \u_t \in \uad_t(\x_t)
  \eqfinp
\label{eq:constraints}
\end{equation}
Note that we do not enforce any explicit bounds on the inner temperature.
Instead, we choose to add a penalization term in the objective function if the temperature
is below a given threshold, as explained below.


\subsection{Objective function}
\label{sec:objective}

At time step~$t$, the instantaneous cost~$L_t: \XX_t \times \UU_t \times \WW_\post \rightarrow \RR$
aggregates two different costs as in the formula
\begin{equation}
  \label{eq:operational_cost}
  L_t(x_t, u_t, w_\post) = \price^e_t  \times \max\na{ 0,f^{ne}_\post } +
  \price^d_t\times  \max\na{0, \overline{\theta^i_t} - \ti_t}
  \eqfinv
\end{equation}
where \( f^{ne}_\post \) is a function of \( \np{u_t, w_\post} \)
by~\eqref{eq:loadbalancerandom},
and \( \ti_t \) is part of the state~$x_t$ in~\eqref{eq:state}.

First, one pays a unitary price~$\price_t^e$ to import electricity from the external network
between step~$t$
and~$t+1$; hence, the electricity cost is equal
to~$\price^e_t \times \max\na{ 0,f^{ne}_\post }  $.
Second, if the indoor temperature is below a given threshold, we penalize the
induced discomfort with a cost
$\price^d_t \times \max\na{0, \overline{\theta^i_t} - \ti_t}$,
where~$\price_t^d$ is a virtual price of discomfort (we choose not
to penalize the temperature in case it is above a given threshold, as we do not
consider any air conditioning in this study).
The cost~$L_t$ is a convex piecewise linear function, a property that will prove important
for the SDDP algorithm.
To ensure that stocks are not empty at the final time step~$\final$,
we add a convex piecewise linear final cost~$K : \XX_{\final} \rightarrow \RR$ of the form
\begin{equation}
  \label{eq:finalpenal}
  K(x_\final) = \kappa \times \max\na{0, x_0 - x_\final}
  \eqfinv
\end{equation}
where~$\kappa$ is a positive penalization coefficient (calibrated by trial and error).

As decisions $\u_t$ and states $\x_t$ are random, the costs $L_t(\x_t, \u_t, \w_\post)$
and $ K(\x_\final)$ become also random variables.
We choose to minimize the expected value of the daily operational cost,
that is,
\begin{equation}
  \label{eq:objective}
  \EE \Bgc{\sum_{t=0}^{\final - 1} L_t(\x_t, \u_t, \w_\post)
  + K(\x_\final) } \eqfinv
\end{equation}
yielding the expected value of a convex piecewise linear cost.

\subsection{Stochastic optimal control formulation}
\label{sec:model:soc}
Finally, the EMS problem is written as a stochastic optimal control 
problem
\begin{subequations}
  \label{eq:stochproblem}
  \begin{align}
    \underset{\x, \u}{\min} ~ & \EE \; \Bgc{\sum_{t=0}^{\final-1}
    L_t(\x_t, \u_t, \alea_{t+1}) + K(\x_\final)} \eqfinv \\
    \st\
    & \x_{t+1} = f_t \bp{\x_t, \u_t, \alea_{\next}}
    \eqsepv  \x_0 = x_0 \eqsepv t=0, \ldots, \final-1
\eqfinv \\
    & \u_t \in \uad_t(\x_t)  \eqsepv t=0, \ldots, \final-1
  \eqfinv \\
    & \sigma(\u_t) \subset \cF_t
 \eqsepv t=0, \ldots, \final-1
\eqfinp
  \label{eq:stochproblem_nonanticipativity}
  \end{align}
\end{subequations}
Problem \eqref{eq:stochproblem} expresses that the microgrid manager
aims to minimize the expected value of the costs while satisfying the dynamics,
the control bounds and the nonanticipativity constraints.

\section{Resolution methods}
\label{sec:algo}

The exact resolution of Problem~\eqref{eq:stochproblem} is out of reach in general.
We propose two different algorithms that provide policies\footnote{%
  See \S\ref{sec:model:cons1} on how the algebraic nonanticipativity
  constraint~\eqref{eq:stochproblem_nonanticipativity} can be turned into a more
  practical functional constraint, making it possible to search for solutions
  that are policies.}
$\pi_t: \WW_1 \times \cdots \times \WW_t \rightarrow \UU_t$ that map available
information~$w_1, \ldots, w_t$ at step~$t$ to a decision~$u_t$.

In~\S\ref{MPC}, we start by presenting how to design management policies
with the MPC algorithm. Then, we depict SDDP-based policies in~\S\ref{sec:sddp}.
Both methods use the dynamics~\( f_t \) in~\eqref{eq:dynamics},
the constraints sets~\( \uad_t(\cdot) \) in~\eqref{eq:constraints},
and the cost functions~$L_t$ in~\eqref{eq:operational_cost}
and~$K$ in~\eqref{eq:finalpenal}.

\subsection{Model Predictive Control (MPC)}
\label{MPC}

Be it ordinary or stochastic, MPC is a classical algorithm
that is commonly used to solve stochastic optimization problems.
Regarding MPC, we follow~\cite{bertsekas1995dynamic}.
At time step~$t$, we consider a deterministic forecast~$(\fw_\post, \ldots,
\fw_\final)$ of the future uncertainties~$(\w_\post, \ldots, \w_\final)$
(see \S\ref{annex:forecast} for more details)
and we solve the following deterministic problem,
where the state~$x_t$ is given:
\begin{subequations}
  \label{eq:mpcproblem}
  \begin{align}
    \underset{(u_t, \cdots, u_{\final-1})}{\min} ~ &
\sum_{j=t}^{\final-1} L_j(x_j, u_j, \fw_{j+1}) + K(x_\final)
\eqfinv \\
    \st\ & x_{j+1} = f_j \bp{x_j, u_j, \fw_{j+1}}
\eqsepv j=t, \ldots, \final-1
\eqfinv\\
    & u_j \in \uad_j(x_j)
\eqsepv j=t, \ldots, \final-1
 \eqfinp
  \end{align}
\end{subequations}
Then, we retrieve the optimal decisions~$(u_t\opt, \ldots, u_{\final-1}\opt)$
and only keep the first decision~$u_t\opt$ to control the system
at the beginning of the time interval~\( \ClosedIntervalOpen{t}{t+1} \).
This procedure is then restarted at step~$t+1$.
  Thus, MPC solves an optimization problem at each time step,
  with a time span going from the current time step~$t$ to the final time step~$\final$.
  Then, at the next time~$t+1$, the optimizer updates the scenario $(\fw_\post, \ldots,
  \fw_\final)$ to take into account the latest observation made at~$t+1$.

\subsection{Stochastic Dual Dynamic Programming (SDDP)}

As said in~\S\ref{sec:intro:litterature},
SDDP is an algorithm widely used to optimize energy systems.
The SDDP algorithm provides, at each time step, a value function and
an online control policy.
Whereas the value functions are computed offline (hence with offline data),
online control policies are computed taking into account an online
probability distribution on the next period noise.

\subsubsection*{Dynamic Programming and Bellman principle}

The Dynamic Programming method~\cite{Bellman57} provides solutions of
Problem~\eqref{eq:stochproblem}
as state feedbacks $\pi_t: \XX_t \rightarrow \UU_t$
(these feedbacks are optimal when the noise process is made of independent
random variables).
Dynamic Programming makes use of a sequence of value functions, obtained offline
by setting $V_\final(x_\final) = K(x_\final)$ and by solving backward in time
the recursive functional equations
\begin{multline}
  \label{eq:bellmanequations}
  V_t(x_t) = \min_{u \in \uad_t(x_t)}
  \int_{\WW_\post} \Bc{ L_t(x_t, u, w_\post) + \\
  V_\post\bp{f_t(x_t, u, w_\post)}}
  \mu^{of}_\post(dw_\post)
\eqfinv
\end{multline}
where~$\mu^{of}_\post$ is a (offline) probability distribution on~$\WW_\post$.
Once these functions obtained, we compute a decision at time step~$t$
as a state feedback
\begin{multline}
  \label{eq:bellmanpolicy}
  \pi_t(x_t) \in \argmin_{u \in \uad_t(x_t)}
  \int_{\WW_\post} \Bc{ L_t(x_t, u, w_\post) + \\
  V_\post\bp{f_t(x_t, u, w_\post)}}
  \mu^{on}_\post(dw_\post) \eqfinv
\end{multline}
where $\mu_\post^{on}$ is an online probability distribution on~$\WW_\post$.
This method proves to be optimal when the random variables~$\w_1, \ldots, \w_\final$
are stagewise independent and when $\mu^{on}_t = \mu^{of}_t$ is the probability distribution
of~$\w_t$.

\subsubsection*{Description of Stochastic Dual Dynamic Programming}
\label{sec:sddp}

Dynamic Programming suffers from the well-known curse of
dimensionality~\cite{Bellman57}: its numerical resolution fails
for state dimension typically greater than~4 when value functions are computed
on a numerical grid.
When uncertainties are stagewise independent random variables, costs~$L_t$ and~$K$ are convex and
dynamics~$f_t$ are linear, the SDDP algorithm provides a solution  to
Problem~\eqref{eq:stochproblem} where the value functions are represented by
convex polyhedral functions~\cite{girardeau2014convergence}.
{
  Most achievements of SDDP have been obtained in the linear case --- when all the costs are linear ---
  as each subproblem~\eqref{eq:bellmanequations} can be solved with the simplex algorithm.
  However, the extension to the generic convex case has not been studied extensively.
}

SDDP provides an outer approximation of the value function
$V_t$ in \eqref{eq:bellmanequations}. Provided a bundle
of $k$ supporting hyperplanes $\na{(\lambda^j_t, \beta^j_t)}_{j = 1, \cdots, k}$
(with, for each $j=1, \cdots, k$, $\lambda^j_t \in \RR^4$ and $\beta^j_t \in \RR$),
the corresponding outer approximation~$\underline{V}_t^k$ is given by
\begin{subequations}
  \begin{align}
    \underline{V}_t^k(x_t) &= \min_{\varrho_t \in \RR} \varrho_t \eqfinv\\
                         \st & \bscal{\lambda_t^j}{x_t} + \beta_t^j \leq \varrho_t
    \eqsepv \forall j = 1, \cdots, k \eqfinp
  \end{align}
\end{subequations}
{
  SDDP takes as input an initial value function $\underline{V}_\post^0$ (usually
$-\infty$) and an initial state $x_0$. In the algorithm, we assume that the probability
distributions have finite support $ \na{w_\post^1, \ldots, w_\post^S}$.
The offline distribution~$\mu_\post^{of}$ is then written in the form $\mu_\post^{of} = \sum_{s=1}^S p_\post^s \delta_{w_\post^s}$,
where $\delta_{w_\post^s}$ is the Dirac measure at~$w_\post^s$
and $p_\post^1, \ldots, p_\post^S$ are probability weights.
}
Then, each iteration~$k$ of SDDP encompasses two passes.
\begin{itemize}
\item
  During the \emph{forward pass}, we draw a scenario~$w_1^k, \ldots,
    w_\final^k$ of uncertainties,
    and obtain a state trajectory~$\ba{x_t^k}_{t = 0 \cdots \final}$
    along this scenario as follows.
    Starting from initial state~$x_0$, we compute~$x_\post^k$ from~$x_t^k$ in an iterative
    fashion: i) we obtain a control~$u_t^k$ at time step~$t$,
    using the available $\underline V_\post^k$ function, by
    \begin{subequations}
    \begin{multline}
      \label{eq:sddpforward}
      u_t^k \in \argmin_{u \in \uad_t(x_t)} \sum_{i=1}^S p_\post^i \bc{
      L_t(x_t^k, u, w_\post^i) +
      \underline V_\post^k\bp{f_t(x_t^k, u, w_\post^i)}}
      \eqfinv
    \end{multline}
    and ii), we set~$x_\post^k = f_t(x_t^k, u_t^k, w_\post^k)$ where $f_t$ is
    the piecewise linear dynamics in~\eqref{eq:dynamics}.
  \item
    During the \emph{backward pass}, we update the approximated value functions
    $\ba{\underline V_t^k}_{t=0, \cdots, \final}$ backward in time along the
    trajectory~$\ba{x_t^k}_{t=0, \cdots, \final}$. At time step~$t$,
    we solve the problem
    \begin{multline}
      \label{eq:sddpbackward}
      \varrho_t^{k+1} =  \min_{u \in \uad_t(x_t)} \sum_{i=1}^S p_\post^i \bc{
      L_t(x_t^k, u, w_\post^i) +
      \underline V_\post^{k+1} \bp{f_t(x_t^k, u, w_\post^i)}}
      \eqfinv
    \end{multline}
    \label{eq:sddpfull}
    \end{subequations}
    and we obtain a new cut~$(\lambda_t^{k+1}, \beta_t^{k+1})$ where~$\lambda_t^{k+1}$
    is a subgradient of the optimal cost function~\eqref{eq:sddpbackward}
    evaluated at the point $x_t = x_t^k$,
    {
    \begin{equation*}
      \lambda_t^{k+1} \in \sum_{i=1}^S p_\post^i \bc{
      \partial_x L_t(x_t^k, u_t^{k+1}, w_\post^i) + \\
      \partial_x \underline V_\post^{k+1} \bp{f_t(x_t^k, u_t^{k+1} w_\post^i)}}
      \eqfinv
    \end{equation*}
  }
    and where $\beta_t^{k+1} = \varrho_t^{k+1} - \bscal{\lambda_t^{k+1}}{x_t^{k}}$.
    This new cut makes it possible to update the function~$\underline V_t^{k+1}$
    by the formula
    $\underline{V}_t^{k+1} = \max\na{\underline V_t^k, \bscal{\lambda_t^{k+1}}{.} + \beta_t^{k+1}}$.
\end{itemize}
We use the stopping criterion introduced in \cite{shapiro2011analysis}
to stop SDDP once the gap between the upper and lower bounds is lower
than 0.1 \%.
{
In practice, the subproblems~\eqref{eq:sddpforward} and \eqref{eq:sddpbackward} are linear,
and can be solved efficiently with any linear programming solver.
At iteration $k$, the lower bound is given directly by $\underline{V}_0(x_0)$.
The upper bound is computed by running $N_{sim}$ forward passes of SDDP
on a fixed set of scenarios $w^1, \cdots, w^{N_{sim}}$; the costs associated with
each scenario are then averaged to get the upper bound (the method is thus akin
to a Monte Carlo simulation).
}
We obtain a sequence $\na{\underline V_t}_{t=0, \cdots,
\final}$ of functions, that are lower approximations of the original Bellman functions.

\subsubsection*{Obtaining online controls with SDDP}

We obtain an online policy by means of the following procedure:
\begin{itemize}
  \item
approximated value functions~$\ba{\underline V_t}$ are computed offline
  with the SDDP algorithm (see~\S\ref{sec:sddp}),
  \item
the approximated value functions~$\ba{\underline V_t}$ are then used
   to compute online a decision at any time step~$t$ for any state~$x_t$ as follows.
\end{itemize}
We compute the SDDP policy~$\pi^{sddp}_t$ by
\begin{multline}
  \label{eq:sddppolicy}
  \pi^{sddp}_t(x_t) \in \argmin_{u \in \uad_t(x_t)} \sum_{i=1}^S p_\post^i \bc{
  L_t(x_t, u, w_\post) + \underline V_\post\bp{f_t(x_t, u, w_\post)}}
  \eqfinv
\end{multline}
which corresponds to replacing the value function~$V_\post$ in
Equation~\eqref{eq:bellmanpolicy} with its approximation $\underline V_\post$.
The decision $\pi_t^{sddp}(x_t)$ is used to control the system between steps~$t$
and~$t+1$. Then, we solve Problem~\eqref{eq:sddppolicy} at step~$t+1$.


\subsection{Discussion}

In this section, we have introduced two methods to design management policies,
the first one based on MPC, the second on SDDP.
Both methods differ on how they model the future uncertainties.
In SDDP, one represents the future as a sequence of independent random variables,
whereas in MPC it is with a deterministic forecast.
It remains now to compare the policy~$\pi^{sdpp}$ with the
policy~$\pi^{mpc}$.

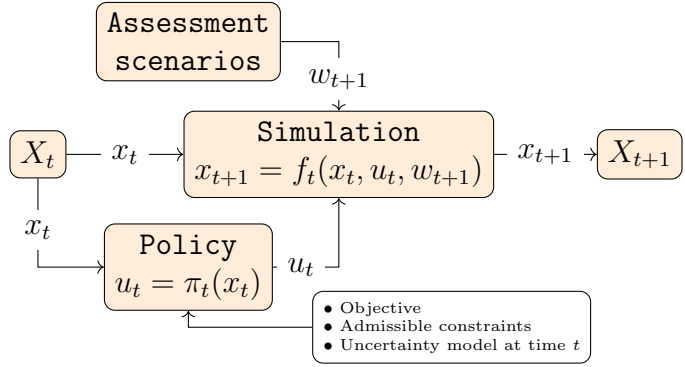
\begin{figure}[!ht]
  \centering

\begin{tikzpicture}[node distance=1.5cm, every node/.style={fill=white, font=\sffamily}, align=center]
  \node (start)             [process]              {Assessment \\ scenarios};
  \node (SimulationBlock)   [process, below of=start,xshift=2cm]
  {Simulation \\ $x_{t+1}=f_t(x_t,u_t,w_{t+1})$};
  \node (XtBlock)           [process, left of=SimulationBlock,xshift=-2.5cm]   {$X_t$};
  \node (Xtp1Block)         [process, right of=SimulationBlock,xshift=+2.5cm]   {$X_{t+1}$};
  \node (PolicyBlock)       [process, below of=XtBlock,xshift=2cm] {Policy \\ $u_t=\pi_t(x_t)$};
  \node (PolicyInfo)  [font=\tiny, align=left, process1,right of=PolicyBlock,xshift=2cm,yshift=-0.8cm]
  {\tiny $\bullet$ Objective \\
    \tiny$\bullet$ Admissible constraints \\
    \tiny$\bullet$ Uncertainty model at time $t$};
  \draw[->]  (PolicyBlock) -| node[xshift=-0.5cm] {$u_t$} (SimulationBlock) ;
  \draw[->]  (XtBlock) |- node[yshift=0.5cm] {$x_t$} (PolicyBlock);
  \draw[->]  (XtBlock) -- node {$x_t$} (SimulationBlock);
  \draw[->] (SimulationBlock) -- node {$x_{t+1}$} (Xtp1Block);
  \draw[->]     (start)  -| node[yshift=-0.5cm] {$w_{t+1}$} (SimulationBlock);
  \draw[->] (PolicyInfo) -| (PolicyBlock) ;
\end{tikzpicture}
  \caption{Flow chart of the simulation procedure}
  \label{fig:simulation_procedure}
\end{figure}

In the next Sect.~\ref{sec:numeric}, we compare the performances
of $\pi^{mpc}$ and $\pi^{sddp}$ on a set of assessment scenarios, using
a simulator. We depict in Figure~\ref{fig:simulation_procedure} the flow chart of the simulation procedure.
The devised policies are used to compute a decision~$u_t$
at each time step~$t$, using the information already available at step~$t$.
Then, by comparing the total costs obtained and repeating
the procedure on a bundle of scenarios, we are able to draw conclusions
about the respective performances of each policy.


\section{Numerical results}
\label{sec:numeric}

In~\S\ref{subsec:casestudy}, we describe a case study.
In~\S\ref{subsec:realistic}, we develop a protocol to fairly
compare the MPC and SDDP algorithms altogether with a rule based policy,
then discuss the results obtained.
In~\S\ref{subsec:whitenoise}, we quantify the robustness
of MPC and SDDP with respect to the level of uncertainty.

\subsection{Case study}
\label{subsec:casestudy}
\paragraph{Settings}
We aim to solve the stochastic optimization problem~\eqref{eq:stochproblem}
over one day, with 96 time steps. The battery size is~$3~\text{kWh}$, and the
hot water tank has a capacity of~$120~\text{l}$.
We suppose that the house has a surface~$A_p = 20~\text{m}^2$ of solar panel
at disposal, oriented south, and with a yield of
15\%. We penalize the recourse variable~$\transfert^{ne}_{t+1}$
in~\eqref{eq:operational_cost} with on-peak and off-peak tariff, corresponding
to \'Electricit\'{e} de France's (EDF) individual tariffs.
The house's thermal envelope corresponds
to the French RT2012 specifications~\cite{rt2012}.
Meteorological data comes from M\'{e}t\'{e}o France measurements corresponding
to the year 2015.

\paragraph{Implementing the algorithms}
We implement MPC and SDDP in Julia 0.6, using JuMP~\cite{dunning2017jump}
as a modeler, \verb+StochDynamicProgramming.jl+ as a SDDP solver, and
Gurobi 7.02~\cite{gurobi2015gurobi} as a LP solver.
All computations run on a Core i7 2.5 GHz processor, with 16GB RAM.

\paragraph{Rule based method}

We choose to compare the MPC and SDDP algorithms with the following basic
decision rule: the battery is charged
whenever the solar production $\PV$ is available, and discharged to fulfill
the demand if there remains enough energy in the battery;
the tank is charged ($\va f^{w}_t > 0$) if the tank energy $\va h_t$ is lower than $\va h_0$,
the heater $\va f^{h}_t$ is switched on when the temperature is below the setpoint $\overline{\theta^i_t}$,
and switched off whenever the temperature is above the setpoint plus a given margin.

\subsection{Benchmark}
\label{subsec:realistic}

\subsubsection*{Demand scenarios}
\begin{figure}[!ht]
  \begin{tabular}{cc}
    \includegraphics[width=5.9cm]{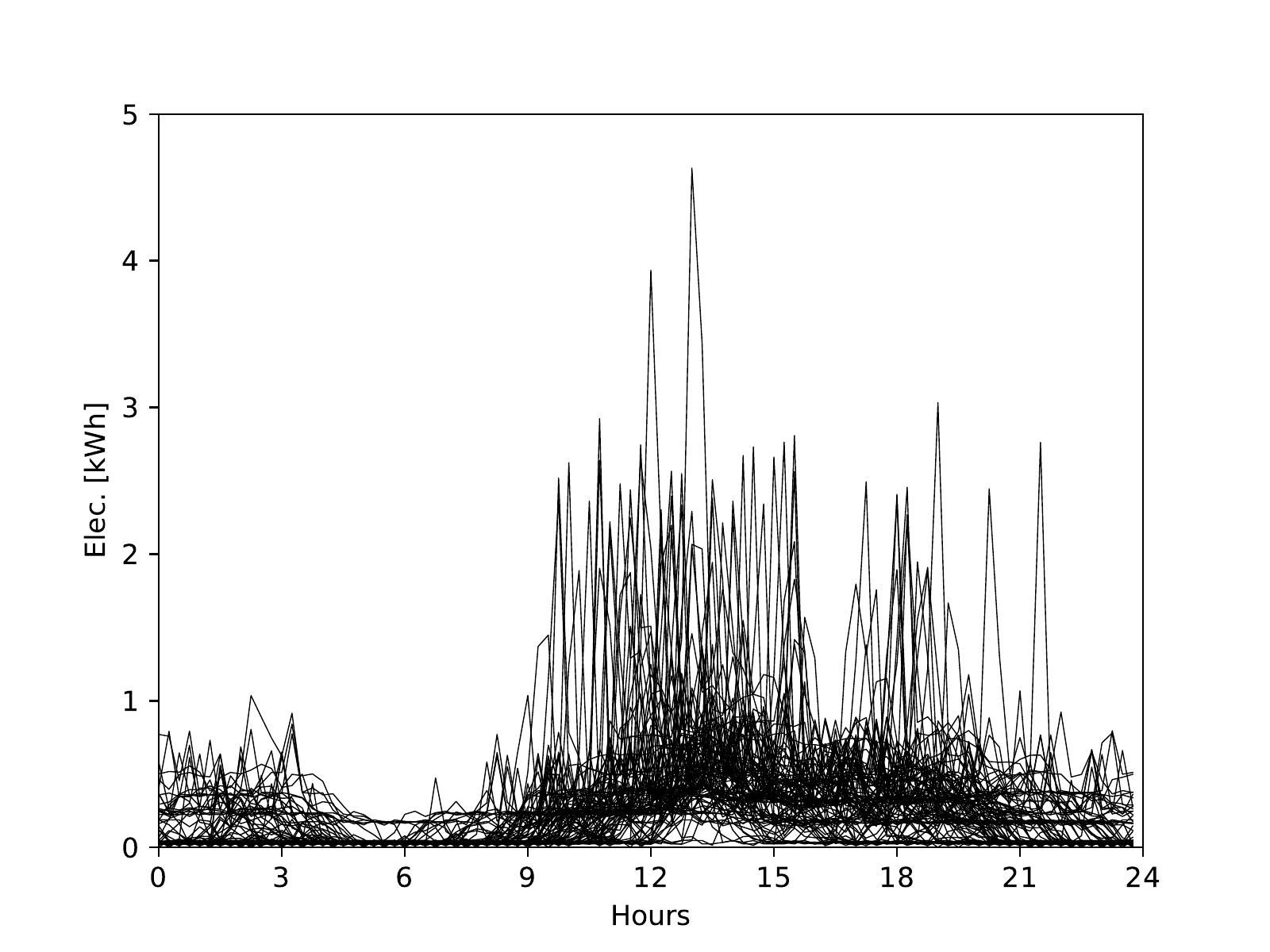} &
    \includegraphics[width=5.9cm]{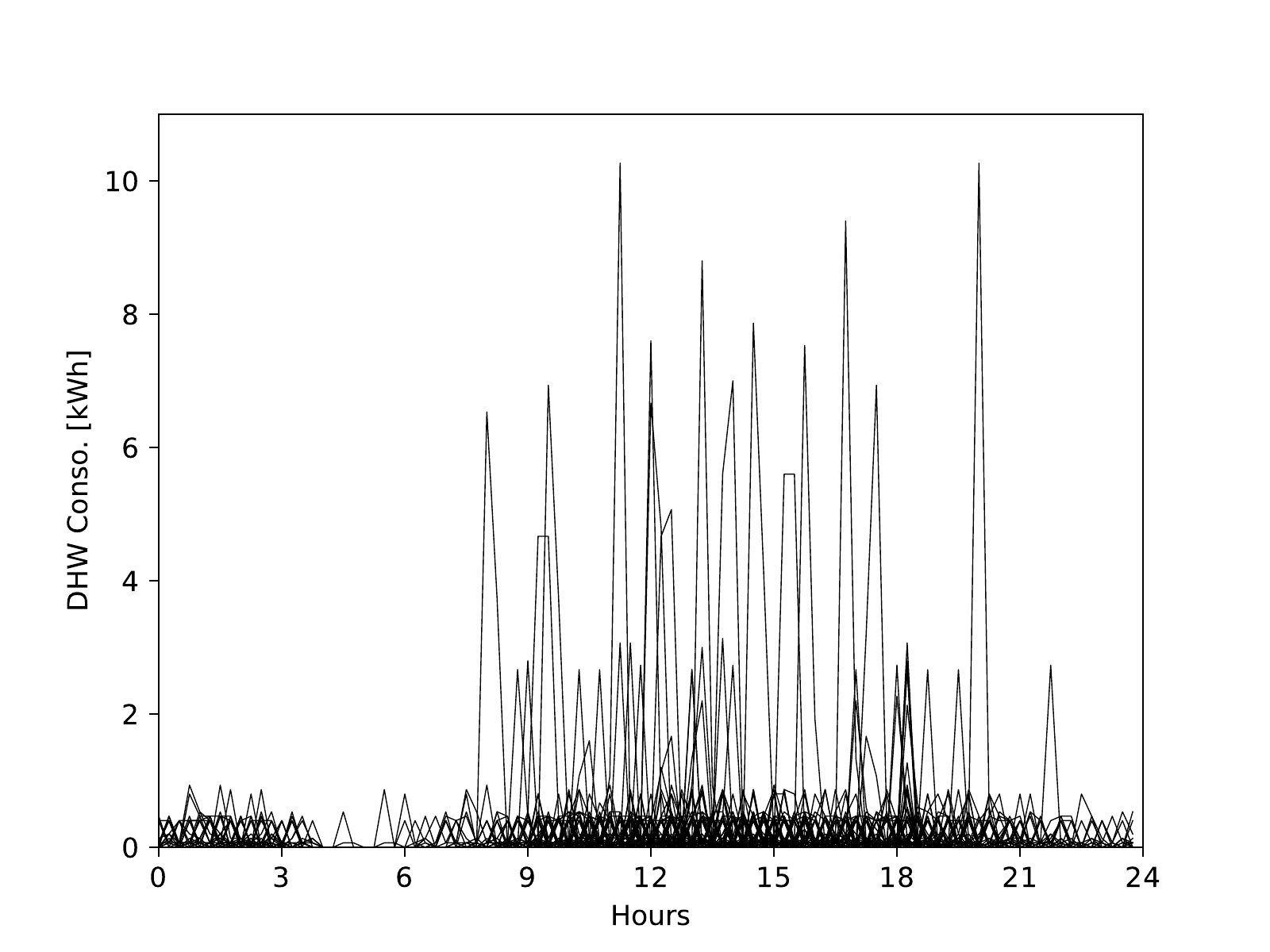}
  \end{tabular}
  \caption{Electrical (left) and domestic hot water (right) demand scenarios}
  \label{fig:scenstrobe}
\end{figure}

Scenarios of electrical and domestic hot water demands, at time steps evenly
spaced every~$\Delta=15$ minutes,
are generated with the software StRoBe~\cite{baetens2016modelling}.
In Figure~\ref{fig:scenstrobe}, 100 scenarios of electrical and hot
water demands are displayed.
We observe almost null demand during the night, and demand peaks around midday and 8~pm;
peaks in hot water demand correspond to showers.
We aggregate the electrical demand~$\Demandel$ minus
the production~$\PV$ of the solar panel
in a single variable~$\Demandel$, so that we consider
only two uncertainties~$\w_t = (\Demandel_t, \va d^{th}_t)$.

{
\subsubsection*{Building offline probability distributions for SDDP}
We use the optimization scenarios to build marginal probability
distributions~$ \mu_t^{of}$ that will feed the SDDP procedure
in~\eqref{eq:sddpforward}-\eqref{eq:sddpbackward}.
In order to obtain a discrete probability distribution at each time~$t$,
we use a Lloyd-Max quantization scheme~\cite{lloyd1982least} to compute
$S$ representative points from the $N_{opt}$ optimization scenarios.
}

\subsubsection*{Out of sample assessment of policies}
To obtain a fair comparison between SDDP and MPC, we use an out-of-sample
validation.
We generate 2,000 scenarios of electrical and hot water demands,
and we split these scenarios in two distinct parts: the first~$N_{opt}= 1,000$
scenarios are called \emph{optimization scenarios}, and the
remaining~$N_{sim}= 1,000$ scenarios are called \emph{assessment scenarios}.

First, during the offline phase, we use the
optimization scenarios to build models for the uncertainties,
under the mathematical form required by each algorithm
(see Sect.~\ref{sec:algo}).
Second, during the online phase, we use the assessment scenarios to compare
the policies produced by these algorithms.
At time step~$t$ during the assessment, the algorithms cannot use the future
values of the assessment scenarios, but can take advantage
of the observed values up to~$t$ to update their statistical models of future uncertainties.
{Our method differs from}
\cite{bhattacharya2018managing,papavasiliou2018application}
where scenarios serve to fit a probability distribution which is used, on the
one hand, to design the SDDP algorithm and, on the other hand, to generate new
scenarios to assess the performances of SDDP.

\subsubsection*{MPC procedure}
Electrical and thermal demands are naturally correlated in
time~\cite{widen2010high}.
To take into account such a dependence across the different
time steps, we choose to model the process~$\w_1,
\ldots, \w_\final$ with an auto-regressive (AR) model and
use the deterministic trend to yield the forecast in~\eqref{eq:mpcproblem}.
We detail the overall procedure in~\S\ref{annex:forecast}.

\subsubsection*{SDDP procedure}

\paragraph{Computing value functions offline}
We fit a sequence of finite probability distributions~$\na{\mu_t^{of}}_{t=0, \cdots, \final}$
using only the optimization scenarios. Then we compute a set of value functions
with the procedure described in~\S\ref{sec:sddp}.

\paragraph{Using value functions online}
Once the value functions have been obtained by SDDP, we
compute online decisions at step~$t$ with Equation~\eqref{eq:sddppolicy},
using a finite online probability distribution~$\mu_t^{on}$, fitted with both
the optimization scenarios and the past of the current assessment scenario
(respecting thus the nonanticipativity constraint).

\subsubsection*{Assessing on different meteorological conditions}
We assess the algorithms on three different days, with different meteorological
conditions (see Table~\ref{tab:meteoday}). Therefore,
we use three distinct sets of $N_{opt} + N_{sim}$ scenarios for demands
(with $N_{opt}$ the number of optimization scenarios and $N_{sim}$ the number
of assessment scenarios), one for each typical day.

\begin{table}[!ht]
 \centering
 \begin{tabular}{lccc}
              & Date           & Temp. ($^\circ C$) & PV Production (kWh) \\
   \hline
   Winter Day & February, 19th & 3.3                & 8.4                 \\
   Spring Day & April, 1st     & 10.1               & 14.8                \\
   Summer Day & May, 31st      & 14.1               & 23.3
  \end{tabular}
  \caption{Different meteorological conditions}
  \label{tab:meteoday}
\end{table}
These three different days correspond to different heating needs. During
Winter day, the heating is maximal, whereas it is medium during Spring day and
null during Summer day. The production of the solar panel varies accordingly.

\subsubsection*{Comparing the algorithms performances}
During assessment, we use MPC (see~\eqref{eq:mpcproblem}) and
SDDP (see~\eqref{eq:sddppolicy}) policies to compute online decisions along
$N_{sim}$ assessment scenarios. Then, we compare the average electricity bill obtained
with these two policies and with the rule based policy.
The assessment results are given in Table~\ref{tab:comparison}, where
means and standard deviation~$\sigma$ are computed with the $N_{sim}=1,000$ 
scenarios; the notation $\pm$ corresponds to the 95\% confidence interval $\pm 1.96
\dfrac{\sigma}{\sqrt{N_{sim}}}$.

\begin{table}[!ht]
  \centering
  \begin{tabular}{lccc}
    \toprule
                      & SDDP            & MPC             & rule based policy \\
    \midrule
    Offline time      & 50~s            & 0~s             & 0~s          \\
    Online time       & 1.5~ms          & 0.5~ms          & 0.005~ms          \\
    \midrule
    Electricity bill (\euro) &                 &                 &           \\
    \midrule
    Winter day        & 4.38~$\pm$ 0.02 & 4.59~$\pm$ 0.02 & 5.55~$\pm$ 0.02          \\
    Spring day        & 1.46~$\pm$ 0.01 & 1.45~$\pm$ 0.01 & 2.83~$\pm$ 0.01          \\
    Summer day        & 0.10~$\pm$ 0.01 & 0.18~$\pm$ 0.01 & 0.33~$\pm$ 0.02          \\
  \end{tabular}
  \caption{Comparison of SDDP, MPC and rule based policies}
  \label{tab:comparison}
\end{table}

Regarding mean electricity bills, we observe, on the three last rows in Table~\ref{tab:comparison},
that MPC and SDDP yield much better results than the rule based policy.
On Winter day, SDDP is slightly better than MPC;
on Spring day, they are equivalent;
on Summer day, SDDP is able to almost halve the mean costs yielded by MPC
(in both case, these costs are low).

\begin{figure}[!ht]
  \centering
  \includegraphics[width=10cm]{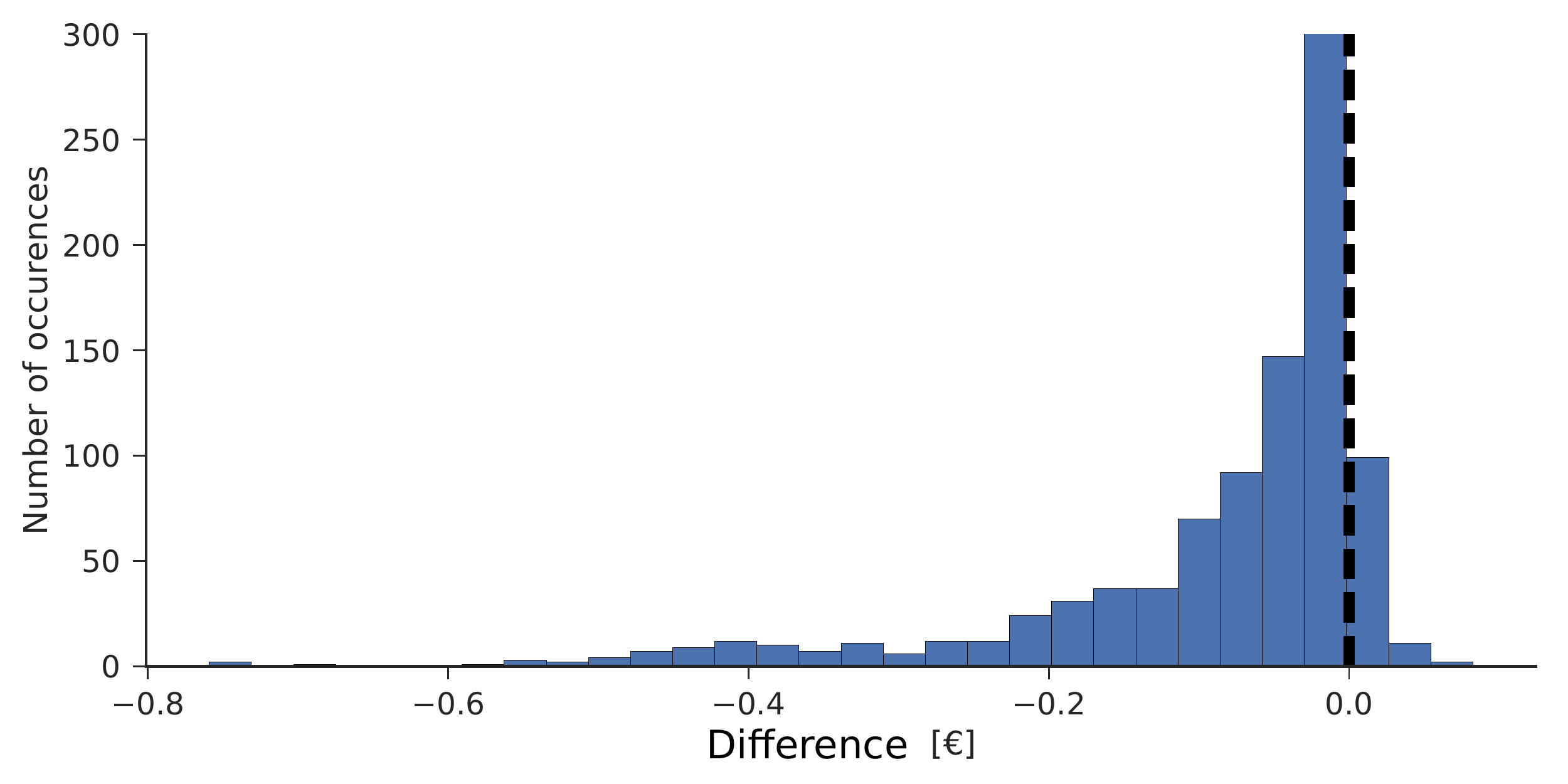}
  \caption{
SDDP total costs minus MPC total costs during Summer day}
  \label{fig:histsavings}
\end{figure}

Figure~\ref{fig:histsavings} displays the histogram of the difference
between SDDP and MPC total costs during Summer day.
{In this way, for each scenario (and not only in the mean), we can measure the discrepancy
between SDDP and MPC performances.}
We observe that SDDP is better than MPC for about 93\% of the scenarios.
The distribution in Figure~\ref{fig:histsavings} exhibits a heavy tail that
reveals the superiority of SDDP on extreme scenarios.
Thus, not only SDDP achieves better performance than MPC in the mean, but also
for the vast majority of scenarios.
Similar analyses hold for Winter and Spring days.



\subsubsection*{Comparing the algorithms running times}

In terms of numerical performance, it takes less than one minute to compute
approximated Bellman functions~$\underline V_t$
as in~\S\ref{sec:sddp} with SDDP on a particular day. Then, the online computation
of a single decision takes 1.5~ms on average, compared to 0.5~ms for MPC.
Indeed, MPC is favored
by the linearity of the optimization Problem~\eqref{eq:mpcproblem},
whereas, for SDDP, the higher the quantization size~$S$ of
the online probability distribution~$\mu_t^{on}$, the slower is the online resolution
of Problem~\eqref{eq:sddppolicy}, but the more accurate $\mu_t^{on}$ is.
MPC's offline resolution time is equal to~0s, as MPC is not an algorithm based
on cost-to-go functions, and as we do not include the offline time
devoted to find coefficients of an AR process.
%
As we are not in a distributed setting, the computation
  time is not a challenge here, and computing the value functions or the
  online decisions is a doable task for a central controller.

\subsubsection*{Analyzing the use of storage capacities}

We analyze now the trajectories of stocks in assessment during Summer day,
where heating is off and production of the solar panel is nominal at midday.

\begin{figure}[!ht]
  \centering
  \includegraphics[width=10cm]{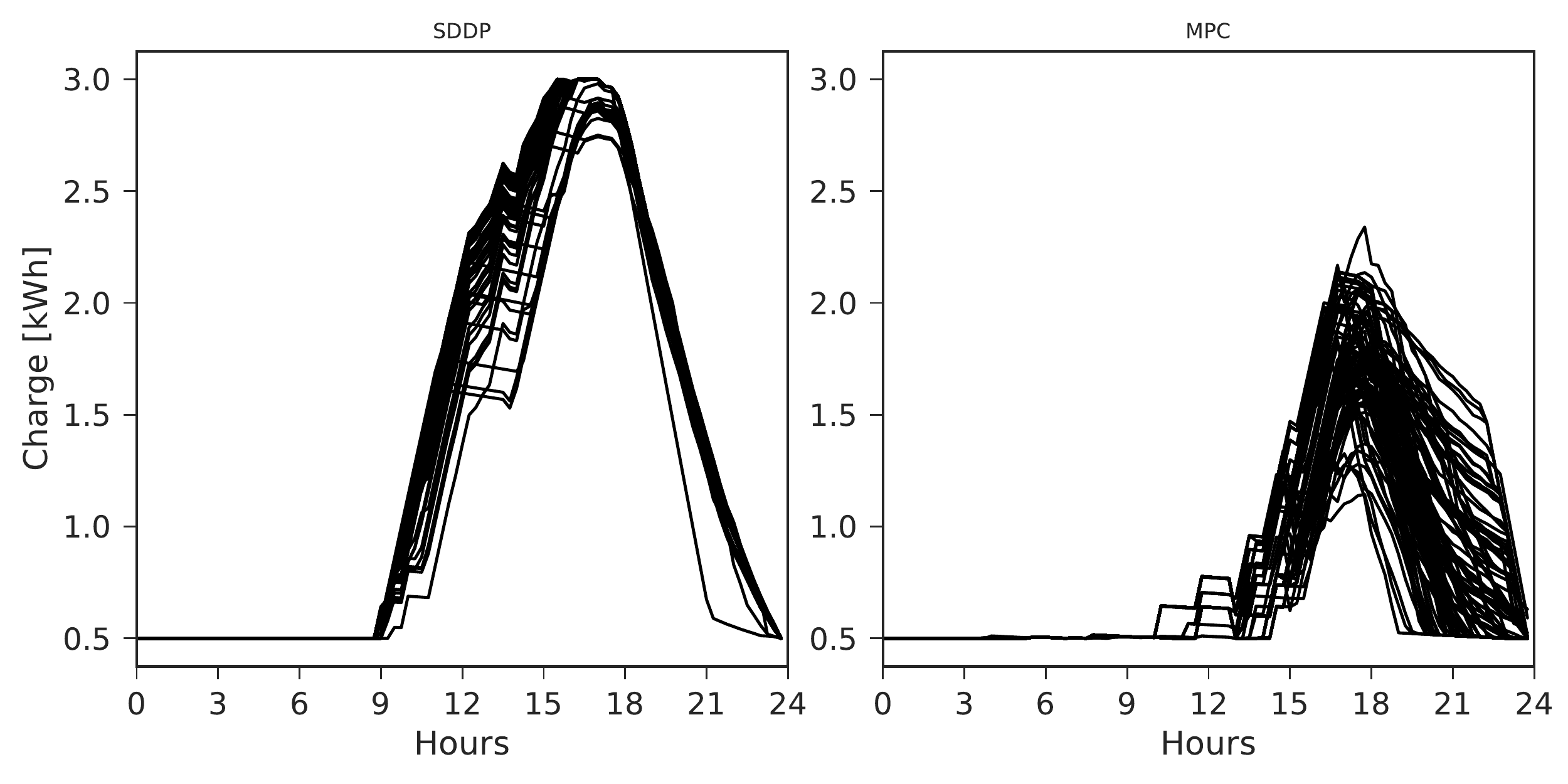}
  \caption{Battery charge trajectories for SDDP (left)
    and MPC (right) during Summer day}
  \label{fig:battraj}
\end{figure}

Figure~\ref{fig:battraj} displays the state of charge of the battery along
a subset of assessment scenarios, for SDDP and MPC. We observe that SDDP
charges earlier the battery at its maximum.
On the contrary, MPC charges
the battery later and does not use the full potential of the battery\footnote{%
In Figure~\ref{fig:battraj}, we observe that SDDP charges the battery to its maximum level during day-time,
when the solar panels are producing electricity. Then, the battery is fully discharged
during the evening to satisfy the local demand when the energy prices are high (evening peak).
On the contrary, MPC does not charge the battery to its maximum level. Indeed, its forecast
only comprises the average energy demand; the algorithm is filling the battery only
to satisfy this average demand (in a sense, we are "overfitting" the average scenario).
Hence, MPC is unable to anticipate a demand higher than usual, leading to an increase use of the recourse variable (the importation from the external grid)
when MPC encounters an unexpected scenario. SDDP does not fall into this pitfall, as it considers
a broader range of demands at each time step with its quantized probability laws (and not just the average demand). This pushes SDDP
to fill the battery to its maximum level, in anticipation of potential high demands during the evening.
  }. The
two algorithms discharge the battery to fulfill the evening demands.
We notice that each trajectory exhibits a single cycle of charge/discharge,
thus preserving the battery's aging.

\begin{figure}[!ht]
  \centering
  \includegraphics[width=10cm]{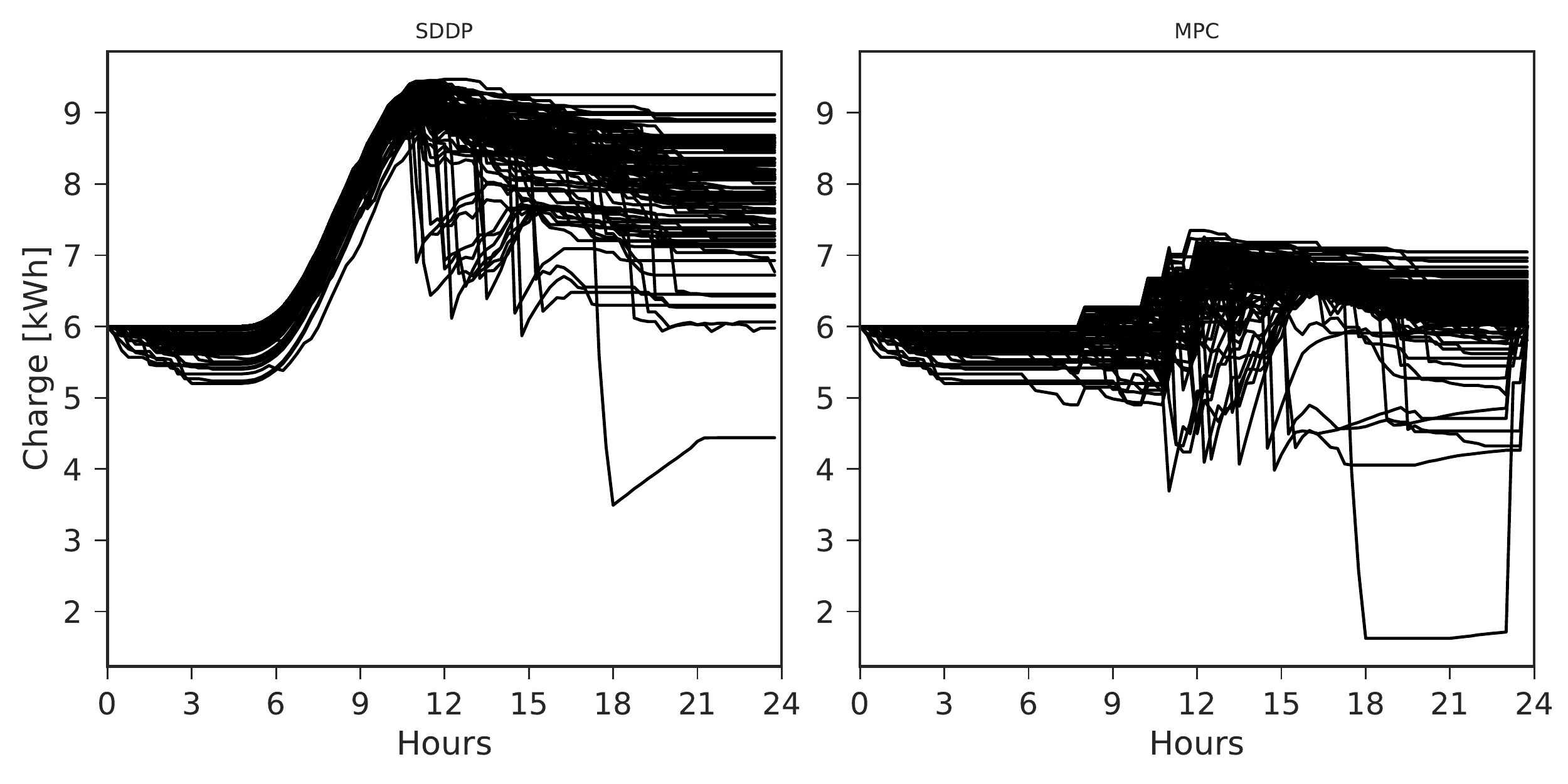}
  \caption{Hot water tank trajectories for SDDP (left) and MPC (right) during Summer day}
  \label{fig:dhwtraj}
\end{figure}

Figure~\ref{fig:dhwtraj} displays the charge of the domestic hot water tank along
the same subset of assessment scenarios. We observe
a similar behavior as for the battery trajectories:
SDDP uses more the electrical hot water tank to store the excess of PV energy,
and the level of the tank is greater at the end of the day than in MPC.

This analysis suggests that SDDP makes a better use of storage capacities than
MPC.

\subsection{Quantifying the sensitivity of MPC and SDDP with respect to the
level of uncertainty}
\label{subsec:whitenoise}


Here, we compare the performances of the two algorithms
(MPC and SDDP) when the level of uncertainty increases.

\paragraph{Uncertainty model}

\begin{figure}[h!]
  \centering
  \begin{tabular}{cc}
    \includegraphics[width=6cm]{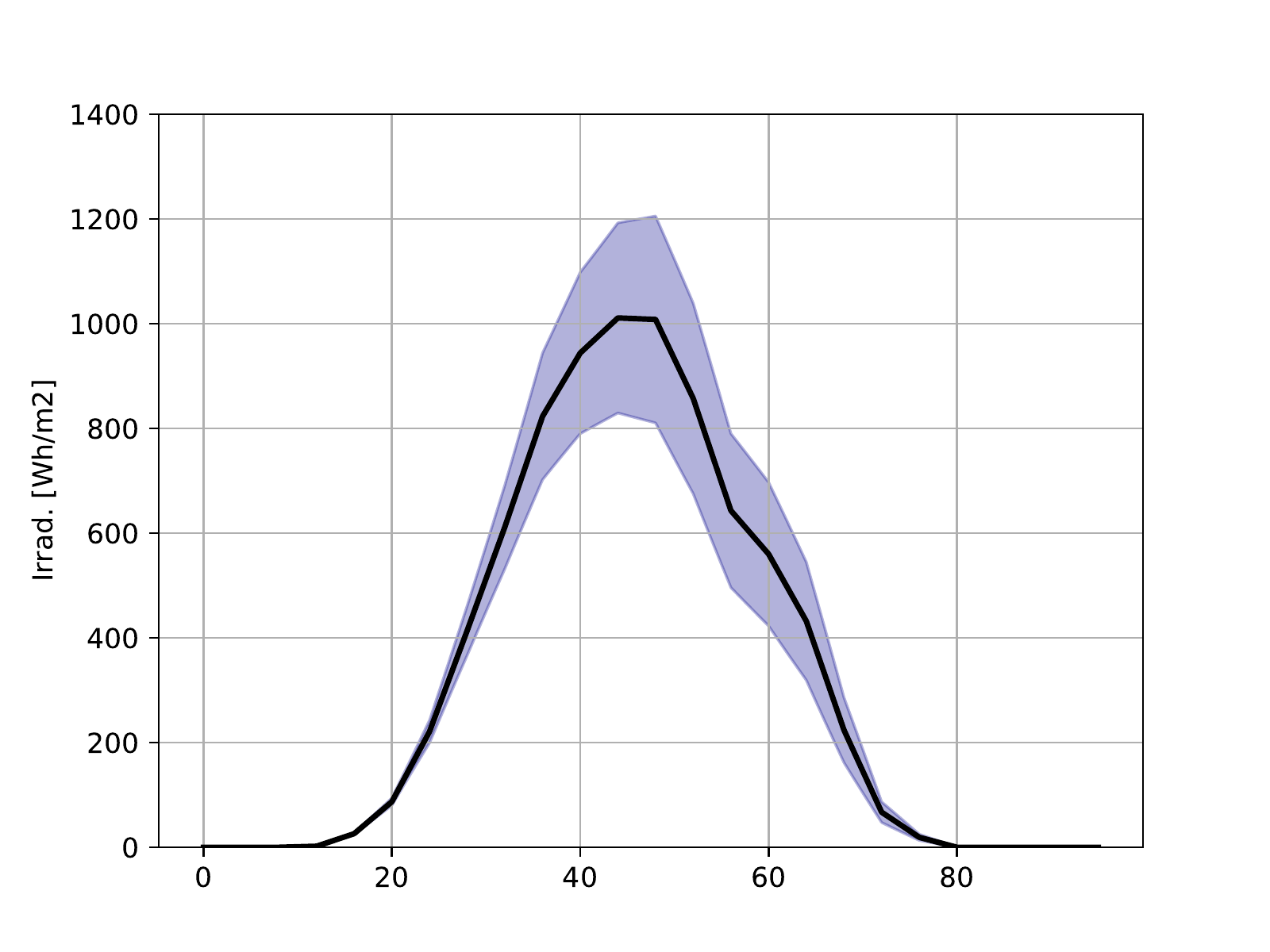}   &
    \includegraphics[width=6cm]{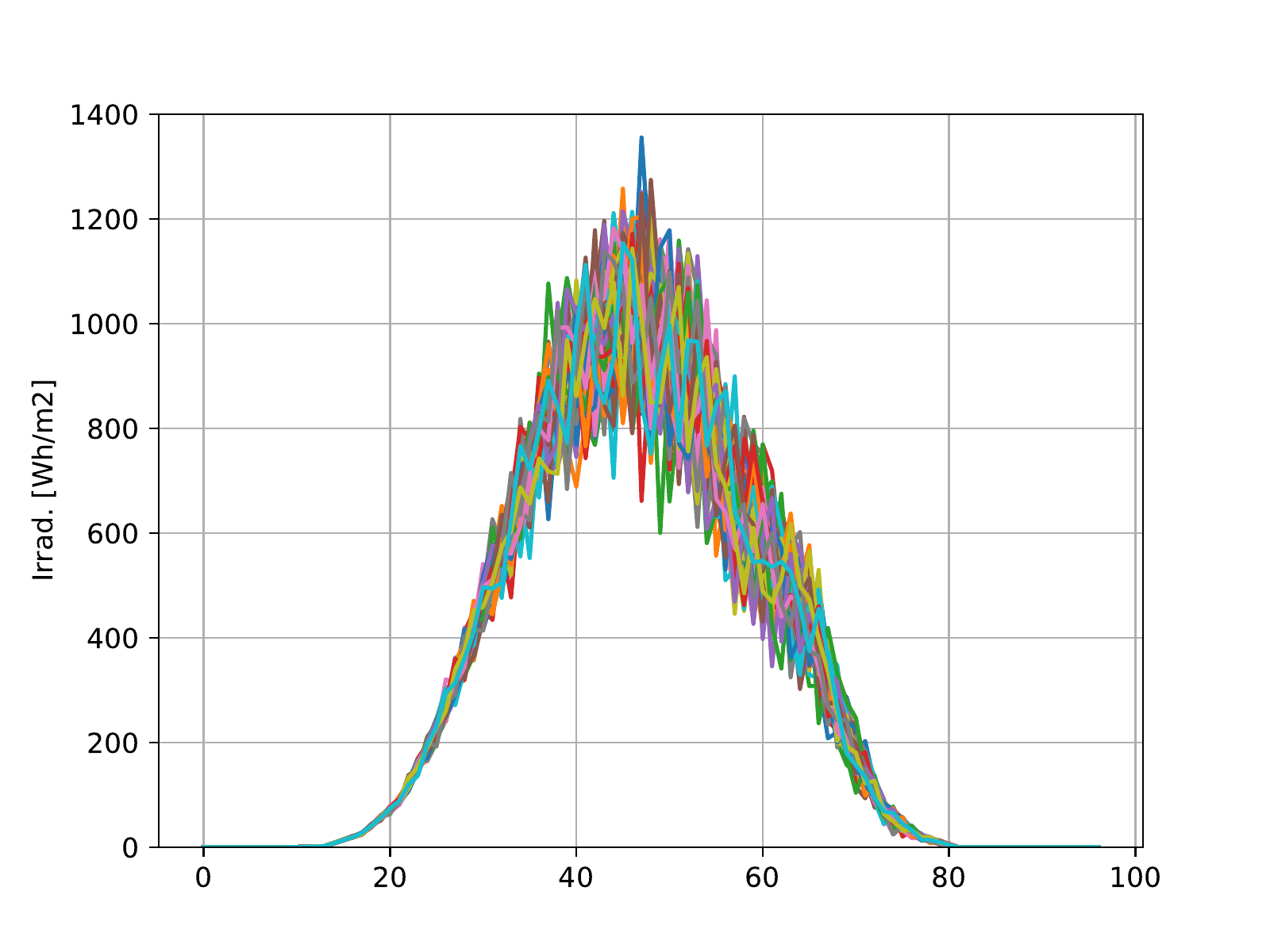} \\
    (a) & (b)
  \end{tabular}
  \caption{Confidence interval (a) and generation of
    corresponding assessment scenarios (b) for
    production of solar panel $\PV$. These scenarios are generated with
    $\sigma_0=0$ and $\sigma_\final = 0.2$.}
    \label{fig:solar:genscen}
\end{figure}

We suppose that, at each time step~$t$, the random variable~$\PV_t$ (the photovoltaic
production in Equation~\eqref{uncertainties}) is given by
\(  \PV_t = \mu_t \times (1 + \va \varepsilon_t) \),
where $\mu_t$ is deterministic (it corresponds to the forecasted value of $\PV_t$)
and where $\va\varepsilon_t$ is a zero-mean Gaussian random variable,
 with standard-deviation
$ \sigma_t = \sigma_0 + (\sigma_\final - \sigma_0) {t}/{\final}$ 
increasing linearly over time
(where $\sigma_0$ is the initial standard-deviation and $\sigma_\final$ the
final standard-deviation). We define the level of uncertainty through the value
of $\sigma_T$: the greater it is, the more difficult it becomes to predict the value of the random variable~$\PV_t$.
With this setting, we have that
$\EE\np{\PV_t} = \mu_t$ and $\vari(\PV_t) = \mu_t^2\sigma_t^2$.
In Figure~\ref{fig:solar:genscen}, we display scenarios that are realizations of
sequences \( \PV_0, \ldots, \PV_\final \).

\paragraph{Results}

\begin{figure}[!ht]
  \centering
    \includegraphics[width=8.5cm]{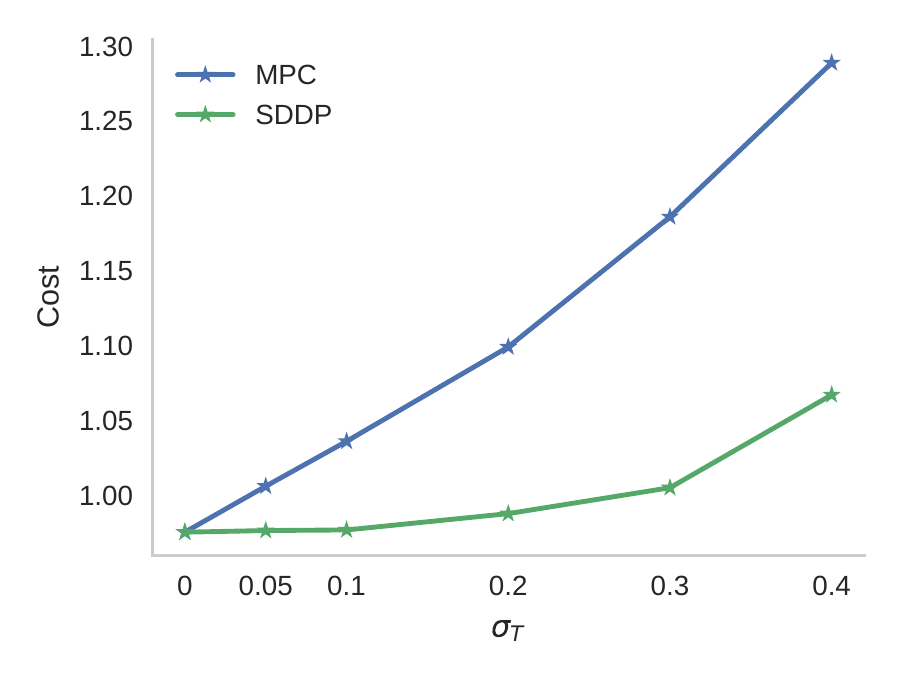}
    \caption{Evolution of the performance of MPC and SDDP as a function of the
level~$\sigma_T$ of uncertainty}
  \label{fig:irradmpcsddp}
\end{figure}

We assess SDDP and MPC with different level of uncertainties,
that is, with increasing values of $\sigma_\final$. The costs correspond
to the management costs to operate the microgrid during a particular
day in Summer, where the production of the solar panel is nominal. The
detailed results are given in Figure~\ref{fig:irradmpcsddp}, which
shows the evolution of the performance of the
two algorithms as a function of the level~$\sigma_T$ of uncertainty.
The costs of each algorithm are obtained via Monte-Carlo simulation, with 10,000
assessment scenarios.
We observe that MPC's cost increases quicker than SDDP's cost
when the level of uncertainty increases.

This second use case allows to quantify the sensitivity
of the algorithms with respect to the level of uncertainty.
It seems that SDDP behaves better than MPC when facing high uncertainties.

\section{Conclusion}
\label{Conclusion}

We have presented the optimal management of a domestic microgrid,
and have compared different management policies
(the core of an Energy Management System) under uncertainties.
Our results show that the two optimization-based policies (MPC and SDDP) outperform the proposed
rule based policy in terms of money savings.
Furthermore, SDDP outperforms MPC during Winter and Summer days ---
and displays similar performance as MPC
during Spring day. Even when SDDP and MPC exhibit close average performances,
a comparison scenario by scenario shows that SDDP beats MPC most of the time
(more than 90\% of scenarios during Summer day).
Thus, SDDP proves better than MPC to manage uncertainties in our study,
although MPC displays also good performances.
SDDP makes better use of storage capacities too.

What if we had compared SDDP not only with ``ordinary'' MPC
but with its stochastic variant, Stochastic Model Predictive Control (SMPC)?
A systematic analysis done in~\cite{LeFranc-Carpentier-Chancelier-DeLara:2020},
using a large dataset of microgrids, reveals that
our conclusions remain valid: algorithms based on the offline computation of
cost-to-go functions (SDP, SDDP) outperform lookahead algorithms
(MPC, SMPC).

Our study can be extended in different directions.
First, we could mix SDDP and MPC to grasp the benefits of these two algorithms.
Indeed, SDDP is designed to handle the uncertainties variability but fails
to capture the time correlation (as it relies on an assumption of stagewise
independence), whereas ordinary MPC ignores the uncertainties variability,
but considers time correlation by means of a future scenario.
Second, we have extended this study in~\cite{Carpentier-Chancelier-DeLara-Pacaud:2020},
where we applied decomposition methods to optimize microgrids comprising
several buildings connected together.
Finally, a comparison of MPC and SDDP with novel methods based on
reinforcement learning could also be of interest.

\section{Appendix}
\label{sec:problem}

In this Appendix, we depict physical equations of the energy system
in Figure~\ref{fig:twostocksschema}. These equations are
naturally written in continuous time~$t$.
We model the battery and the hot water tank with stock dynamics, and
the dynamics of the house's temperatures with an electrical analogy.



\subsection{Energy storage}
\label{sec:twostockselec}

We consider a battery, whose state of charge at time~$t$ is
denoted by~$b(t)$.
The battery dynamics is given by the differential equation
\begin{equation}
  \label{eq:batteryequation}
  \dfrac{db}{dt} = \rho_c (f^b(t))^+ - \dfrac{1}{\rho_d} (f^b(t))^- \eqfinv
\end{equation}
with~$\rho_c$ and~$\rho_d$ being the charge and discharge efficiency
and~$f^b(t)$ denoting the energy exchange with the battery.


\subsection{Electrical hot water tank}

We use a simple linear model for the electrical hot water tank dynamics.
The enthalpy balance equation writes
\begin{equation}
    \label{eq:tankequation}
    \dfrac{dh}{dt} = \beta_\tank f^{w}(t) - d^{th}(t) \eqfinv
\end{equation}
where
\begin{itemize}
  \item $f^{h}(t)$ is the electrical energy used to heat the tank, satisfying
    \begin{equation}
    0 \leq f^{w}(t) \leq \overline f^{w} \eqfinv
    \end{equation}
  \item $d^{th}(t)$ is the domestic hot water demand,
  \item $\beta_\tank$ is a conversion yield.
\end{itemize}

\subsection{Thermal envelope}
We model the evolution of the temperatures inside the house with
an electrical analogy: we view temperatures as voltages, walls as capacitors,
and thermal flows as currents. A model with 6 resistances and 2 capacitors
(R6C2) proves to be accurate to describe small buildings \cite{berthou2014development}.
The model takes into account two temperatures:
\begin{itemize}
  \item the wall's temperature~$\tw(t)$,
  \item the inner temperature~$\ti(t)$.
\end{itemize}
Their evolution is governed by the two following differential
equations
\label{sec:twostockstherm}
\begin{subequations}
  \label{eq:r6c2continuous}
  \begin{multline}
    c_m \dfrac{d \tw}{dt} = \underbrace{\dfrac{\ti(t) - \tw(t)}{R_i + R_s}}_{\substack{\text{Exchange} \\ \text{  Indoor/Wall}}}
    + \underbrace{\dfrac{\to(t) - \tw(t)}{R_m + R_e}}_{\substack{\text{Exchange} \\ \text{  Outdoor/Wall}}}
    + \underbrace{\gamma f^{t}(t)}_{\text{Heater}}\\
    + \underbrace{\dfrac{R_i}{R_i + R_s} \pint(t)}_{\substack{\text{Radiation} \\ \text{through windows}}}
    + \underbrace{\dfrac{R_e}{R_e + R_m} \pext(t)}_{\substack{\text{Radiation} \\ \text{through wall}}} \eqfinv
  \end{multline}
  \begin{multline}
    c_i \dfrac{d \ti}{dt} = \underbrace{\dfrac{\tw(t) - \ti(t)}{R_i + R_s}}_{\substack{\text{Exchange} \\ \text{  Indoor/Wall}}}
    + \underbrace{\dfrac{\to(t) - \ti(t)}{R_v}}_{\text{Ventilation}}
    + \underbrace{\dfrac{\to(t) - \ti(t)}{R_f}}_{\text{Windows}}\\
    + \underbrace{(1 - \gamma) f^{t}(t)}_{\text{Heater}}
    + \underbrace{\dfrac{R_s}{R_i + R_s} \pint(t)}_{\substack{\text{Radiation} \\ \text{through windows}}} \eqfinv
  \end{multline}
\end{subequations}
where we denote
\begin{itemize}
  \item the energy injected in the heater by~$f^{h}(t)$,
  \item the external temperature by~$\to(t)$,
  \item the radiation through the wall by~$\pext(t)$,
  \item the radiation through the windows by~$\pint(t)$.
\end{itemize}
The time-varying quantities~$\to(t)$, $\pint(t)$ and~$\pext(t)$ are exogenous.
We denote by~$R_i, R_s, R_m, R_e, R_v, R_f$ the different resistances of the
R6C2 model,
and by~$c_i, c_m$ the capacities of the inner rooms
and the walls. We denote by $\gamma$ the proportion of heating dissipated in the wall
through conduction, and by $(1- \gamma)$ the proportion of heating dissipated
in the inner room through convection.
We detail the numerical values in Table~\ref{tab:elecanalogy}.

\begin{table}[!ht]
  \centering
  \begin{tabular}{cc}
    \hline
    $R_i$ & $4.81 \times 10^{-4}$~SI \\
    $R_s$ & $2.94 \times 10^{-4}$~SI \\
    $R_m$ & $4.51 \times 10^{-3}$~SI \\
    $R_e$ & $1.48 \times 10^{-4}$~SI \\
    $R_v$ & $4.51 \times 10^{-3}$~SI \\
    $R_f$ & $2.00 \times 10^{-2}$~SI \\
    $c_i$ & $8.30 \times 10^7$ ~SI \\
    $c_m$ & $5.85 \times 10^6$ ~SI \\
    \hline
  \end{tabular}
    \caption{Numerical values for the electrical analogy\label{tab:elecanalogy}}
\end{table}

\subsection{MPC}
\label{annex:forecast}

\paragraph{Building offline an AR model for MPC}

We fit an AR(1) model using the optimization scenarios
(we do not consider higher order lags for the sake of simplicity).
For $i \in \{el, hw\}$, the AR model writes
\begin{subequations}
\begin{equation}
  d^{i}_{t+1} = \alpha_t^{i} d^{i}_{t} + \beta_t^{i} + \varepsilon_t^{i} \eqfinv
\end{equation}
where the nonstationary coefficients~$(\alpha^i_t, \beta^i_t)$ are, for any time step~$t$,
solutions of the least-square problem
\begin{equation}
  (\alpha^i_t , \beta_t^i) = \argmin_{a,b}\sum_{s=1}^{N_{opt}}
  \norm{d^{i,s}_\post - ad^{i, s}_{t} - b}^2_2 \eqfinp
\end{equation}
The points $(d^{i, 1}_t, \ldots
d_t^{i, N_{opt}})$ correspond to the optimization scenarios.
The AR residuals $(\varepsilon_t^{el}, \varepsilon_t^{hw})$ are a white noise process.
\end{subequations}

\paragraph{Updating the forecast online}
Once the AR model is calibrated, we use it to update the forecast during
assessment (see~\S\ref{MPC}). The update procedure is threefold:
\begin{itemize}
  \item[i)] we observe the demands~$w_t = (d^{el}_t, d^{hw}_t)$ between time steps~$t-1$ and~$t$,
  \item[ii)] we update the forecast~$\overline w_{t+1}$ at time step~$t+1$ with the AR model
    \begin{equation*}
      \overline w_\post =
      \bp{\overline d^{el}_\post, \overline d^{hw}_\post} =
      \Bp{\alpha_t^{el} d^{el}_{t} + \beta_t^{el}, \;
      \alpha_t^{hw} d^{hw}_{t} + \beta_t^{hw} }
      \eqfinv
    \end{equation*}
  \item[iii)] we set the forecast between time steps~$t+2$ and~$\final$ by
    using the mean values
    \begin{equation*}
      \overline w_\tau = \dfrac{1}{N_{opt}} \sum_{i=1}^{N_{opt}} w_\tau^{i}
      \quad \forall \tau = t+2, \cdots, \final
    \end{equation*}
 of the optimization scenarios
\end{itemize}
Once the forecast $(\fw_\post, \ldots, \fw_\final)$ is available,
it serves as input into the optimization Problem~\eqref{eq:mpcproblem}
(the MPC algorithm).

\end{document}